\documentclass[12pt, oneside]{article}   	% use "amsart" instead of "article" for AMSLaTeX format
\usepackage{geometry}                		% See geometry.pdf to learn the layout options. There are lots.
\usepackage{graphicx}				% Use pdf, png, jpg, or eps§ with pdflatex; use eps in DVI mode
\usepackage{amssymb}
\usepackage{cite}
\usepackage{latexsym}
\usepackage{amsmath}

\def\rr{\mathbb{R}^{N}}
\def\ep{\varepsilon}
\def\n{\rho}
\def\de{\delta}
\def\ga{\gamma}
\def\na{\nabla}
\def\p{\partial}
\def\vp{\varphi }
\def\th{\theta}

\newtheorem{theorem}{Theorem}[section]
\newtheorem{proposition}{Proposition}[section]
\newtheorem{remark}{Remark}[section]
\newtheorem{definition}{Definition}[section]
\newtheorem{lemma}[theorem]{Lemma}

\newcommand{\bn}{\begin{eqnarray}}
\newcommand{\en}{\end{eqnarray}}
\newcommand{\ba}{\begin{aligned}}
\newcommand{\ea}{\end{aligned}}
\newcommand{\bnn}{\begin{eqnarray*}}
\newcommand{\enn}{\end{eqnarray*}}
\newcommand{\be}{\begin{equation}}
\newcommand{\ee}{\end{equation}}
\newcommand{\la}{\label}

\title{Vanishing pressure  limit     for   compressible Navier-Stokes equations with degenerate viscosities}

\author{Zhilei L{\small IANG}\thanks{The author   is
  supported  by  NNSFC Grant No.11301422. E-mail:  zhilei0592@gmail.com}  \\[3mm]
{\normalsize    School of Economic Mathematics,} \\
{\normalsize Southwestern  University of Finance and Economics, Chengdu  611130,  China} \\[2mm]}

 \date{}	 % Activate to display a given date or no date

\begin{document}
\maketitle

\begin{abstract}
In this paper we study   a  vanishing pressure process for highly  compressible Navier-Stokes  equations as the Mach number tends to infinity. We first prove  the  global existence of  weak solutions for the pressureless  system   in the framework  [Li-Xin, arXiv:1504.06826v2], where  the  weak solutions   are established for   compressible Navier-Stokes  equations with degenerate viscous  coefficients.  Furthermore,  a  rate  of convergence of   the  density in $L^{\infty}\left(0,T;L^{2}(\rr)\right)$ is obtained,  in   case when   the  velocity  corresponds  to the gradient of  density at initial time.
\end{abstract}

\section{Introduction}
 The time evolution of a  viscous compressible barotropic  fluid occupying the whole space $\rr\, (N=2,\,3)$ is governed by the  equations
\be\ba\la{1}
\p_{t}\rho_{\ep} +{\rm div}(\rho_{\ep} u_{\ep}) = 0,\\
\p_{t}(\rho_{\ep}u_{\ep}) +  {\rm div}(\n_{\ep}u_{\ep} \otimes u_{\ep}) + \na P_{\ep}-{\rm div \mathbb{S}_{\ep}}=0,\ea\ee where the unknown functions $\rho_{\ep}$ and  $u_{\ep}$ are the density and the velocity. The pressure $P_{\ep}= \ep\n_{\ep}^{\gamma} $ with $\gamma>1$ is given and  $\ep>0$ is related to Mach number, and the stress tensor takes the form
\be\la{e6}  \mathbb{S}_{\ep}= h(\n_{\ep})\na u_{\ep} + g(\n_{\ep}) {\rm div} u_{\ep}\mathbb{I},\ee
in which $\mathbb{I}$ is the identical matrix, $ h$ and $g$ are functions of $\n_{\ep}$ satisfying  the physical restrictions \be\la{e6e} h(\n_{\ep})>0,\quad h(\n_{\ep})+Ng(\n_{\ep})\ge0.\ee
For   simplicity reason, in this paper we assume
\be\la{6} h(\n_{\ep})=\n_{\ep}^{\alpha},\quad g(\n_{\ep})=(\alpha-1)\n_{\ep}^{\alpha}, \quad{\rm for}\,\,\,\alpha> (N-1)/N,\quad \ep\in (0,1).\ee
The initial functions are imposed as
\be\la{5} \n_{\ep}(x,t=0)=\n_{0}\ge0,\quad \n_{\ep} u_{\ep}(x,t=0)=m_{0},\quad x\in \rr.\ee

For fixed $\ep,$   equations \eqref{1} is one of the most important mathematical models describing the motion of a viscous flow.   There is a huge literature     on the existence and asymptotics  of solutions  because of its   mathematical  challenges  and  wide physical  applications; see  \cite{desjar1,enp,mv,yu,kim,lixin,llx,p2} and the references therein.  Consider the constant viscosities  $h$ and $g$ defined in \eqref{e6e}, and more general symmetric stress tensor $$\mathbb{S}_{\ep}= h \left(\frac{\na u_{\ep}+(\na u_{\ep})^{tr}}{2}\right)+ g  {\rm div} u_{\ep}\mathbb{I},$$   Lions\cite{p2} first proved the global  existence of weak solutions of \eqref{1} if the  adiabatic index $\ga\ge 3N/(N+2)$. Later,  $\ga$ was relaxed by Feireisl-Novotny-Petzeltov$\acute{a}$ \cite{enp}  to $\ga>N/2$ and by Jiang-Zhang \cite{jiangzhang}  to $\ga>1$ under some  extra spherically symmetry assumptions.  While for the  case when  viscosities are  density-dependent,   Bresch-Desjardins \cite{bres,bres1}   developed
a new entropy structure on  condition that
\bnn g(\n_{\ep})=\n_{\ep} h^{'}(\n_{\ep})-h(\n_{\ep}).\enn
This  gives  an  estimate on the gradient of  density, and thereby,  some further compactness information on density.
Li-Li-Xin \cite{llx} proposed the global entropy weak solution to system \eqref{1}  in one-dimensional bounded interval  and studied  the vacuum vanishing   phenomena  in  finite time span.  Similar results in \cite{llx} were extended to the Cauchy problem in \cite{jiu} by Jiu-Xin.  Guo-Jiu-Xin \cite{guo} obtained the global existence of weak solution to \eqref{1} if  some   spherically symmetric assumptions are made.     However, the problem becomes much more difficult in  general  high dimension spaces.
 Mellet-Vasseur \cite{mv} provided a compactness framework which ensures the existence of weak solutions as a  limit of  approximation solutions, but leaves  such   approximations    sequence open  in \cite{mv}.  Until recently, the problem was  solved in two impressing  papers by Vasseur-Yu \cite{yu}  and Li-Xin   \cite{lixin},   where they   constructed separately  appropriate approximations  from  different approaches.  Vasseur-Yu \cite{yu1}  also  considered the compressible quantum Navier-Stokes equations with damping, which helps  to understand  the existence of global weak solutions to the compressible Navier-Stokes equations.

Refer to \cite{mv,lixin}, we give  the weak solution of system \eqref{1} in  below
  \begin{definition} \la{def1} For   fixed  $\ep>0,$  we call    $(\n_{\ep},u_{\ep})$  a weak solution to the problem  \eqref{1}-\eqref{5}, if
\bnn\begin{cases}
0\le \rho_{\ep}\in L^{\infty}\left(0,T; L^{1}(\rr)\cap L^{\ga}(\rr)\right),\\
\na \n_{\ep}^{\alpha-1/2},\,\, \sqrt{\n_{\ep}}u_{\ep}\in L^{\infty}\left(0,T; L^{2}(\rr)^{N}\right),\\
\na \n_{\ep}^{(\ga+\alpha-1)/2}\in L^{2}\left(0,T; L^{2}(\rr)^{N}\right),\\
h(\n_{\ep}) \na u_{\ep}\in L^{2}\left(0,T; W_{loc}^{-1,1}(\rr)^{N\times N}\right),\\
g(\n_{\ep}) {\rm div}u_{\ep}\in L^{2}\left(0,T; W_{loc}^{-1,1}(\rr)\right),\\
\end{cases} \enn
 $(\sqrt{\n_{\ep}}, u_{\ep})$
 satisfy   $\eqref{1}_{1}$ in distribution sense, and   the  integral equality
 \bnn \ba \int_{\rr}m_{0}\phi(x,0)+&\int_{0}^{T}\int_{\rr}\sqrt{\n_{\ep}}(\sqrt{\n_{\ep}} u_{\ep})\p_{t}\phi+\sqrt{\n_{\ep}}u_{\ep} \otimes \sqrt{\n_{\ep}}u_{\ep}: \na \phi+\ep \n_{\ep}^{\ga}{\rm div}\phi\\
 =&<h(\n_{\ep})\na u_{\ep}, \na \phi>+<g(\n_{\ep}){\rm div}u_{\ep},{\rm div}\phi>\ea\enn
holds  true for any   test functions  $\phi\in C_{0}^{\infty}\left(\rr\times [0,T)\right)^{N}$,
where
 \bnn\ba <h(\n_{\ep})\na u_{\ep}, \na \phi>=-\int_{0}^{T}\int_{\rr}\left(\n_{\ep}^{\alpha-1/2}\sqrt{\n_{\ep}}u_{\ep} \triangle \phi+\frac{2\alpha}{2\alpha-1}\sqrt{\n_{\ep}}u_{j\ep}\p_{i}\n_{\ep}^{\alpha-1/2}\p_{i}\phi_{j}\right)\ea\enn
 and
 \bnn\ba& <g(\n_{\ep}){\rm div}u_{\ep},{\rm div}\phi>\\
 &=-(\alpha-1)\int_{0}^{T}\int_{\rr}\left(\n_{\ep}^{\alpha-1/2}\sqrt{\n_{\ep}}u_{\ep} \cdot \na {\rm div} \phi
 +\frac{2\alpha}{2\alpha-1} \sqrt{\n_{\ep}}u_{\ep}\na \n^{\alpha-1/2}{\rm div}\phi\right).\ea\enn
 \end{definition}

The  following  important  existence results of weak  solutions  are obtained in   \cite{lixin} by  Li-Xin.
 \begin{proposition}(\cite{lixin})
\la{pro}
Assume that the initial function in \eqref{5} satisfy
\be\la{e}\begin{cases} 0\le \n_{0}\in L^{1}(\rr)\cap L^{\ga}(\rr),\quad \na \n_{0}^{\alpha-1/2}\in L^{2}(\rr),\\
m_{0}\in L^{2\ga/(\ga+1)}(\rr),\quad  \n_{0} \not\equiv 0,\,\,\,m_{0}=0\,\,a.e. \,\,{\rm on}\,\, \{x\in \rr|\, \n_{0}=0\},\\
\n_{0}^{-(1+\eta_{0})}|m_{0}|^{2+\eta_{0}}\in L^{1}(\rr)\quad {\rm for\,\,some}\,\,\,\eta_{0}>0.\end{cases}\ee
Additionally,   assume  that  for $N=2$
\be\la{e1} \alpha>1/2,\,\,\,\gamma>1, \,\,\,\,\gamma \ge2\alpha-1,\ee
and for  $N=3$
\be\la{e2} \begin{cases} \ga\in (1,3),\\
  \gamma\in (1,6\alpha-3) \quad {\rm if} \,\,\alpha \in [3/4,1],\\
 \gamma\in  [2\alpha-1,  3\alpha-1]\,\,\,{\rm and}\,\,\,  \rho_{0}^{-3}|m_{0}|^{4} \in L^{1}(\mathbb{R}^{3}) \quad{\rm if}\,\,  \alpha\in (1,2).\end{cases} \ee

 Then   the problem  \eqref{1}-\eqref{5} has global weak solutions $(\rho_{\ep},u_{\ep})$ in the sense of Definition \ref{def1}.\end{proposition}
 \begin{remark} \la{r0} By   \eqref{e2}, if we  multiply   $\eqref{1}_{2}$ by $4|u_{\ep}|^{2}u_{\ep}$  and compute  directly, we infer
\be\la{e4} \n_{\ep}|u_{\ep}|^{4}\in L^{\infty}\left(0,T;L^{1}(\mathbb{R}^{3})\right).\ee
{\em Proof.}  The rigorous  proof is available in  Appendix.
 \end{remark}

It seems
rather natural to expect that, as  $\ep\rightarrow 0+$, the limit   $(\n,u)$  of $(\rho_{\ep},u_{\ep})$ satisfy   the corresponding  pressureless system
 \be\ba\la{1a}
\p_{t}\rho+{\rm div}(\rho u) = 0,\\
\p_{t}(\rho u) +  {\rm div}(\n u \otimes u)-{\rm div}\left(h(\n)\na u+ g(\n) {\rm div} u\mathbb{I}\right)=0.\ea\ee As in  \cite{has1}, we define the weak solution (called {\em quasi-solution}) to system \eqref{1a}
\begin{definition} \la{def2} The function   $(\n,u)$ is   called a {\em quasi-solution}, if
\bnn\begin{cases}
0\le \rho\in L^{\infty}\left(0,T; L^{1}(\rr)\right),\\
\na \n^{\alpha-1/2},\,\, \sqrt{\n}u\in L^{\infty}\left(0,T; L^{2}(\rr)^{N}\right),\\
h(\n) \na u\in L^{2}\left(0,T; W_{loc}^{-1,1}(\rr)^{N\times N}\right),\\
g(\n) {\rm div}u\in L^{2}\left(0,T; W_{loc}^{-1,1}(\rr)\right);\\
\end{cases} \enn
 in addition,  $(\sqrt{\n}, u)$
 satisfy   $\eqref{1a}_{1}$ in distribution sense, and the  integral equality
 \bnn \ba \int_{\rr}m_{0}\phi(x,0)+&\int_{0}^{T}\int_{\rr}\sqrt{\n}(\sqrt{\n} u)\p_{t}\phi+\sqrt{\n}u \otimes \sqrt{\n}u: \na \phi\\
 =&<h(\n)\na u, \na \phi>+<g(\n){\rm div}u,{\rm div}\phi>\ea\enn
holds  true, where  the quantities on the right  side are     defined as the same of  $<h(\n_{\ep})\na u_{\ep}, \na \phi>$ and $<g(\n_{\ep}){\rm div}u_{\ep},{\rm div}\phi>$. \end{definition}

\begin{remark}The   {\em quasi-solution} in Definition \ref{def2}  was  first proposed  by B. Haspot  to approximate in some sense the compressible Navier-Stokes equations. See, for example, the paper \cite{has1}.\end{remark}

In this   paper, we choose  $\ep=\eta^{-2}>0$ with $\eta$ being  the Mach number. The  readers can refer  to   the papers such as \cite{gou,has1, lions,p2} for more information in  this aspect.  There are   satisfactory results on the incompressible limit when $\eta\rightarrow0,$  we refer readers to  the pioneer works by Desjardins-Grenier-Lions-Masmoudi\cite{desjar,desjar1, lions} when the viscous coefficients are constant. Regretfully,  seldom  result is available up to publication when  $\eta\rightarrow \infty$.  One major difficulty is the  compactness lack of the      density because  its  $L^{\ga}$-bound   is no longer conserved  for  constant viscosities. However,  the case of  density-dependent viscosity  is    much  different  due to the new  BD entropy inequality.   Haspot \cite{has1} proved   the highly compressible limit ($\ep\rightarrow0$) in  the sense of distribution in
suitable Lebesgue spaces, and discussed the  global existence  of   $quasi$-$solutions$  as a convergence limit from  approximation solutions  of system  \eqref{1}, although    such approximations are only \emph{a-priori} exist. It is worthy   mentioning    that  in \cite{has1} the author constructed a family of  explicit solutions $(\n,u)$  with   $\n$   satisfying  the porous medium equation, the heat equation, or the fast diffusion equation,  up to  the choice of $\alpha$.  Haspot-Zatorska \cite{has} consider the one-dimensional Cauchy problem and obtain a rate of convergence of $\n_{\ep}$ and other related properties.

We are interested in  the limit procedure for  the weak solutions $(\n_{\ep}, u_{\ep})$ of  \eqref{1} as $\ep$ tends to zero, and then get some     convergence rate the solutions.  In particular, on the basis of existence results obtained in  \cite{lixin} by Li-Xin,  we adopt some ideas in  \cite{has1} and \cite{mv} and first show the     $quasi$-$solutions$  stability  for the   solutions $(\n_{\ep}, u_{\ep})$ of \eqref{1}.  Secondly,  in
the spirit of  \cite{has}, we obtain  a convergence rate of $\n-\n_{\ep}$ in terms of $\ep$ in high dimensions  by the argument of duality, as long as  the initial velocity  associated    with the gradient of  initial density.

\begin{theorem}
\la{t1} Let    the conditions \eqref{e}-\eqref{e2}  in Proposition \ref{pro} hold true.   Then, for $\alpha\ge 1$,  the  solution  $(\n_{\ep}, u_{\ep})$ of \eqref{1} converges to a limit function  $(\n,u)$  which  solves  \eqref{1a} in the sense of Definition \ref{def2}.
Furthermore,  \be\la{bb1} \n_{\ep}\rightarrow \n \quad   {\rm in}\quad C\left([0,T];L_{loc}^{q_{1}}(\rr)\right),\ee
\be\la{bb2} \n_{\ep}u_{\ep}\rightarrow \n u\quad   {\rm in}\quad L^{2}\left(0,T;L^{q_{2}}_{loc}(\rr)\right),\ee
 where $q_{1}\in [1,\infty),\,\,q_{2}\in [1,2)$ if $N=2$;  $q_{1}\in [1,6\alpha-3),\,\,q_{2}\in [1,\frac{12\alpha-6}{6\alpha-1})$ if $N=3$.\end{theorem}

  \begin{remark} \la{r0} The assumptions in Proposition \ref{pro}  guarantee  the existence of $(\n_{\ep},u_{\ep})$ to \eqref{1}, whose proof are available in \cite{lixin}.  We  allow  more general viscosities at the cost of  stress tensor having the form \eqref{e6}, although it seems   not appropriate from a physical point of view. \end{remark}
 \begin{remark} \la{r1} In case of $\alpha=1,$  Theorem \ref{t1} is valid for  the symmetric viscous stress tensor  $\mathbb{S}_{\ep}={\rm div} \left(\n_{\ep} \frac{\na u_{\ep}+(\na u_{\ep})^{tr}}{2}\right),$ where the existence of $(\n_{\ep},u_{\ep})$ are achieved in \cite{lixin,yu}.  Moreover,   the case  $\alpha<1$ can also be discussed  by  modifying  slightly  the argument in Theorem \ref{t1}.  \end{remark}

\begin{theorem} \la{t2}  In addition to the  assumptions made in Theorem \ref{t1}, let     \be\la{11}u_{0}+\alpha \n_{0}^{\alpha-2}\na \n_{0}=0.\ee   Then there is a positive  $C$ which may depend on  $T$ such that  for $\alpha\ge 3/2$
\bnn  \sup_{0\le t\le T}\|(\n_{\ep}-\n)(\cdot,t)\|_{L^{2}(\rr)}\le C\ep^{\sigma},\enn
 where $\sigma<\frac{1}{2(2\alpha-1)}$ if $N=2$,  $\sigma=\frac{4\alpha-3}{4(2\alpha-1)^{2}}$ if $N=3.$
\end{theorem}

 \begin{remark} \la{r3} For one-dimensional problem,  Haspot-Zatorska  \cite{has} first obtained    a   rate of convergence of $\n_{\ep}-\n$ in suitable Sobolev spaces for    $1<\alpha\le 3/2$.  We remark that the argument in \cite{has} relies heavily on the upper bound of density.    \end{remark}

In the rest of this paper, Section 2   is for   some useful lemmas,  and Sections 3-4 are devoted  to  proving  Theorem \ref{t1} and Theorem \ref{t2}  respectively.

\section{Preliminaries}
\begin{lemma}(see \cite{Gar,la}) \la{lem2.2} Let $B_{R} =\{x\in \rr\,:\, |x|<R\}$.  For any $v\in W^{1,q}(B_{R})\cap L^{r}(B_{R}),$ it satisfies
\be\la{lpo}\ba \|v\|_{L^{p}(B_{R})}\le C_{1}\|v\|_{L^{r}(B_{R})}+ C_{2}\|\na v\|_{L^{q}(B_{R})}^{\ga}\|v\|_{L^{r}(B_{R})}^{1-\ga},\ea\ee where   the constant $C_{i} (i=1,2)$  depends only  on $p,q,r, \ga;$   and  the exponents $0\le \ga\le 1,$ $1\le q,\,r\le \infty $  satisfy
$\frac{1}{p}=\ga(\frac{1}{q}-\frac{1}{N})+(1-\ga)\frac{1}{r}$  and
\bnn\left\{
  \begin{array}{ll}
\min\{r,\,\frac{Nq}{N-q}\}\le p\le \max\{r,\,\frac{Nq}{N-q}\}, & {\rm if}\,\,q<N; \\
 r\le p<\infty, &  {\rm if}\,\,q=N;\\
 r\le p\le \infty, &  {\rm if}\,\,q>N.
  \end{array}
\right.\enn
 \end{lemma}

The following $L^{p}$-bound estimate is taken from  \cite[Lemma 2.4]{ll},   whose proof is available by adopting \cite[Lemma 12]{kim} and the  elliptic theory due to Agmon-Douglis-Nirenberg \cite{doug}.
\begin{lemma} (\cite[Lemma 2.4]{ll}) \la{lem2.2} Let $p\in (1,+\infty)$ and  $k\in \mathbb{N}$. Then for  all $v\in W^{2+k,p}(B_{R})$
with  0-Dirichlet boundary condition, it holds that
\bnn \|\na^{2+k}v\|_{L^{p}(B_{R})}\le C\|\triangle v\|_{W^{k,p}(B_{R})},\enn
where the $C$  relies  only on $p$ and $k.$
 \end{lemma}

\begin{lemma}  \la{lem2.4}  Assume  that $f$ is increasing and convex in $\mathbb{R}_{+}=[0,+\infty)$  with  $f(0)=0$. Then,
\bnn |x-y|f(|x-y|)\le (x-y)(f(x)-f(y)),\quad \forall \,\,\,\, x,\,y\in \mathbb{R}_{+}.\enn
\end{lemma}
{\it Proof.} Define $F(x)=f(x)-f(y)-f(x-y).$    Since $f$ is convex, then   $F'(x)\ge 0$ for $x\ge y\ge0$.  This  and  $F(0)=0$ deduce   $F(x)\ge  0$. Hence,
\bnn f(|x-y|)=f(x-y)\le f(x)-f(y)=|f(x)-f(y)|.\enn
Repeating the argument when  $y\ge x\ge0$, we obtain
\bnn f(|x-y|)\le |f(x)-f(y)|,\quad \forall \,\,\,\, x,\,y\in \mathbb{R}_{+}.\enn
This, along with the monotonicity  of $f$, leads to
\bnn\ba   |x-y|f(|x-y|)&\le |x-y| |f(x)-f(y)| = (x-y)(f(x)-f(y)),\ea\enn
the required.

\section{Proof of Theorem \ref{t1}}
In what follows,   the  operations  are based on   hypotheses  imposed   in  Proposition \ref{pro}, and the generic  constant $C>0$ is     $\ep$ independent.

Firstly,   for all existing time $t\ge 0,$ we have
\be\la{70} \|\n_{\ep}(\cdot,t)\|_{L^{1}(\rr)}=\|\n_{0}\|_{L^{1}(\rr)} \ee
and \be\ba\la{7}&\sup_{0\le t\le T}\int_{\rr}\left(\n_{\ep}|u_{\ep}|^{2}+\ep \n_{\ep}^{\ga}\right)+\int^{T}_{0}\int_{\rr}\n_{\ep}^{\alpha}|\na u_{\ep}|^{2}+(\alpha-1)\n_{\ep}^{\alpha}({\rm div} u_{\ep})^{2}\\
&\le \int_{\rr}\frac{|m_{0}|^{2}}{\n_{0}}+\ep \int_{\rr}\n_{0}^{\ga}.\ea\ee
Following  in \cite{mv},  a straight calculation  shows \be\ba\la{8} &\sup_{0\le t\le T}\int_{\rr}\left(\n_{\ep}|u_{\ep}+\alpha \n_{\ep}^{\alpha-2}\na \n_{\ep}|^{2}+\ep \n_{\ep}^{\ga}\right)+\alpha \ga\ep\int^{T}_{0}\int_{\rr}\n_{\ep}^{\alpha+\ga-3}|\na \n_{\ep}|^{2}\\
&\le \int_{\rr}\n_{0}|u_{0}+\alpha \n_{0}^{\alpha-2}\na \n_{0}|^{2}+\ep \int_{\rr}\n_{0}^{\ga}.\ea\ee
The initial condition \eqref{e} and  \eqref{70}-\eqref{8} guarantee
\be\la{15s0}\ba \sup_{0\le t\le T}\int_{\rr}\left(\n_{\ep} + \left|\na \n_{\ep}^{\alpha-1/2}\right|^{2}\right)dx\le C.\ea\ee

We claim that
 \be\ba \la{15s} &\n_{\ep}\in L^{\infty}\left(0,T; L^{1}(\rr)\cap L^{q}(\rr) \right),\quad  \na \n_{\ep}^{\alpha-1/2}\in L^{\infty}\left(0,T; L^{2}(\rr)\right),\ea \ee
where    $q<\infty$ if $N=2,$ and $q=6\alpha-3$ if $N=3$.\\
\underline{{\em Proof of   \eqref{15s}}}. If  $\alpha \le 3/2,$
by \eqref{lpo} we have
\be\la{s02}\ba \| \n_{\ep}^{\alpha-1/2}\|_{L^{p}(\rr)}&\le C\left(\| \n_{\ep}^{\alpha-1/2}\|_{L^{(\alpha-1/2)^{-1}}(\rr)}+ \|\na  \n_{\ep}^{\alpha-1/2} \|_{L^{2}(\rr)} \right)\\
& \le C\left(\| \n_{\ep}\|_{L^{1}(\rr)}^{\alpha-1/2}+ \|\na  \n_{\ep}^{\alpha-1/2} \|_{L^{2}(\rr)} \right),\ea\ee  where $p\ge (\alpha-1/2)^{-1}$ for $N=2$ and $p=6$ for $N=3.$ While for  $\alpha>3/2,$ by \eqref{lpo} and   interpolation theorem, one has
\be\la{s02c}\ba \| \n_{\ep}^{\alpha-1/2}\|_{L^{p}(\rr)}&\le C\left(\| \n_{\ep}^{\alpha-1/2}\|_{L^{1}(\rr)}+C\|\na  \n_{\ep}^{\alpha-1/2} \|_{L^{2}(\rr)} \right)\\
&\le C\left(\| \n_{\ep}\|_{L^{1}(\rr)}^{(1-\th)(\alpha-1/2)}\|\n_{\ep}^{\alpha-1/2}\|_{L^{p}(\rr)}^{\th}+ \|\na  \n_{\ep}^{\alpha-1/2} \|_{L^{2}(\rr)}\right),\ea\ee  where $\th=\frac{p(\alpha-1/2)-p}{p(\alpha-1/2)-1}$,   $p>1$ if  $N=2$ and $p=6$ if $N=3.$
The  \eqref{15s} thus  follows from \eqref{15s0},
  \eqref{s02} and \eqref{s02c}.

The key  issue  in proving  Theorem\ref{t1} is to get the $\ep$-independent estimates and   take $\ep$-limit in Definition \ref{def1}.  In terms of  \eqref{e2} and \eqref{15s}, one has
\bnn \ep\int_{0}^{T}\int_{\rr}\n_{\ep}^{\ga}{\rm div} \phi \rightarrow 0\quad {\rm as}\quad \ep\rightarrow0.\enn
Besides, we also need to   justify  \eqref{bb1}, \eqref{bb2} and  the strong convergence of $\sqrt{\n_{\ep}} u_{\ep}.$    For that purpose it suffices  to   prove  the    Lemmas \ref{lem3.1}-\ref{lem3.2} below.

\begin{lemma} \la{lem3.1} Upon  to some  subsequence,  it satisfies
 \be\la{y35}\n_{\ep}\rightarrow \n \quad   {\rm in}\quad C\left([0,T];L_{loc}^{q_{1}}(\rr)\right),\ee where
 $q_{1}\in [1,\infty)$ if $N=2$ and $q_{1}\in [1,6\alpha-3)$ if $N=3$.      \end{lemma}
\emph{Proof.} By  \eqref{15s} and H$\ddot{o}$lder inequality, we have
\be\la{s5} \|\na \n_{\ep}^{\alpha}\|_{L^{1}(\rr)}\le C\|\n_{\ep}\|_{L^{1}(\rr)}^{1/2}\|\na \n_{\ep}^{\alpha-1/2}\|_{L^{2}(\rr)}\le C.\ee
Since   $1\le \alpha\le 2\alpha-1<6\alpha-3$, from \eqref{15s} and  \eqref{7}
we deduce \bnn\ba &\|\n_{\ep}^{\alpha}u_{\ep}\|_{L^{1}(\rr)}+\int_{0}^{T}\|\n_{\ep}^{\alpha} {\rm div} u_{\ep}\|_{L^{1}(\rr)}\\
&\le \|\n_{\ep}^{\alpha-1/2}\|_{L^{2}(\rr)}\|  \sqrt{\n_{\ep}}u_{\ep}\|_{L^{2}(\rr)}\\
&\quad+\sup_{0\le t\le T}\|\n_{\ep}^{\alpha/2}\|_{L^{2}(\rr)}\int_{0}^{T}\|\n_{\ep}^{\alpha/2} {\rm div} u_{\ep}\|_{L^{2}(\rr)}\le C.\ea\enn This, along  with  \be\la{y13}\p_{t}\n_{\ep}^{\alpha}=(1-\alpha)\n_{\ep}^{\alpha} {\rm div} u_{\ep}-{\rm div} (\n_{\ep}^{\alpha} u_{\ep}), \ee ensures  that $\p_{t}\n_{\ep}^{\alpha}\in L^{2}\left(0,T;W_{loc}^{-1,1}(\rr)\right).$  By the  Aubin-Lions Lemma, we get
\bnn \n_{\ep}^{\alpha}\rightarrow \n^{\alpha} \quad {\rm in}\,\,\,C\left([0,T];L_{loc}^{\beta}(\rr)\right)\quad {\rm for}\,\,\,\beta\in [1,3/2).\enn
Therefore, up to some subsequence,   \be\la{s53} \n_{\ep}^{\alpha}\rightarrow \n^{\alpha},\quad almost\,\,everywhere.\ee
So,  \eqref{15s} and \eqref{s53} guarantee  the strong convergence of $\n_{\ep}$ to $\n$ in  $L^{\infty}\left(0,T;L_{loc}^{\underline{q}}(\rr)\right)$ with
 $\underline{q}\in [1,\infty)$ if $N=2$ and $\underline{q}\in [1,6\alpha-3)$ if $N=3$. Choosing $\alpha=1$ implies $\p_{t}\n_{\ep}\in L^{2}\left(0,T;W_{loc}^{-1,1}(\rr)\right),$ we conclude  \eqref{y35} by   the Aubin-Lions lemma.

  Consequently,  the \eqref{15s} and  \eqref{s53} implies that
   \bnn \sqrt{\n_{\ep}}\rightharpoonup \sqrt{\n}\quad{\rm in}\,\,L^{2}\left(0,T;L^{2}_{loc}(\rr)\right),\quad
 \n_{\ep}^{\alpha-1/2}\rightharpoonup  \n^{\alpha-1/2}\quad{\rm in}\,\,L^{2}\left(0,T;H^{1}_{loc}(\rr)\right).\enn

\begin{lemma} \la{lem3.3}  Upon  to some  subsequence,  it satisfies \be\la{y40} \sqrt{\n_{\ep}}u_{\ep}\rightarrow \sqrt{\n}u\quad{\rm in}\,\,L^{2}\left(0,T;L^{2}_{loc}(\rr)\right).\ee\end{lemma}
\emph{Proof.} The process is divided  into several steps.

\emph{Step 1.}  Define  $m_{\ep}=\left(\chi(\n_{\ep})\n_{\ep}^{\alpha}+(1-\chi(\n_{\ep}))\n_{\ep}^{(1+\alpha)/2}\right)u_{\ep}$ with  $\chi(x)$ being smooth and satisfying   $\chi(x)=1\,\,{\rm}\,\,if\,\, |x|\le 1$ and $ \chi(x)=0\,\,{\rm}\,\,if \,\, |x|\ge2.$

 We claim that, for  $p\in [1,3/2),$
\bnn\ba
& m_{\ep} \rightarrow \left(\chi(\n)\n^{\alpha}+(1-\chi(\n))\n^{(1+\alpha)/2}\right)u\quad{\rm in}\,\, L^{2}\left(0,T;L^{p}_{loc}(\rr)\right). \ea\enn Consequently,
\be\la{y20}
 m_{\ep}\rightarrow \left(\chi(\n)\n^{\alpha}+(1-\chi(\n))\n^{(1+\alpha)/2}\right)u,  \,\,\,almost\,\,everywhere.\ee

In fact, we deduce from  \eqref{15s} and \eqref{7}  that
\bnn\ba &\int_{0}^{T}\|\na(\chi(\n_{\ep})\n_{\ep}^{\alpha}u_{\ep})\|_{L^{1}(\rr)}^{2}\\
&\le C\int_{0}^{T}\left(\|\n_{\ep}^{\alpha}|u_{\ep}||\na \n_{\ep}| \|_{L^{1}(\{1\le \n_{\ep}\le 2\})}^{2}+
 \|\chi\sqrt{\n_{\ep}}|u_{\ep}|\na\n_{\ep}^{\alpha-1/2}| +\chi  \n_{\ep}^{\alpha}|\na u_{\ep}|\|_{L^{1}(\rr)}^{2} \right)\\
&\le C\int_{0}^{T}\left(\|\sqrt{\n_{\ep}}u_{\ep} \|_{L^{2}}^{2}\|\na\n_{\ep}^{\alpha-1/2} \|_{L^{2}(\rr)}^{2}+\|\n_{\ep}^{\alpha/2}\|_{L^{2}(\rr)}^{2}
\|\n_{\ep}^{\alpha/2}\na u_{\ep}\|_{L^{2}(\rr)}^{2}\right)\\
&\le C
\ea\enn and
\be\la{y17}\ba &\int_{0}^{T}\|\na((1-\chi(\n_{\ep}))\n_{\ep}^{(1+\alpha)/2}u_{\ep})\|_{L^{1}(\rr)}^{2}dt\\
&\le C\int_{0}^{T} \|\n_{\ep}^{(1+\alpha)/2}|u_{\ep}||\na \n_{\ep}| \|_{L^{1}(\{1\le \n_{\ep}\le 2\})}^{2}\\
&\quad +C\int_{0}^{T}\|(1-\chi(\n_{\ep}))(\n_{\ep}^{(2-\alpha)/2}|u_{\ep}|\na  \n_{\ep}^{\alpha-1/2}|+\n_{\ep}^{(1+\alpha)/2}|\na u_{\ep}|)\|_{L^{1}(\{1\le \n_{\ep}\})}^{2}\\
&\le C\int_{0}^{T}\left(\|\sqrt{\n_{\ep}}u_{\ep} \|_{L^{2}(\rr)}^{2}\|\na\n_{\ep}^{\alpha-1/2} \|_{L^{2}(\rr)}^{2} +\|\sqrt{\n_{\ep}}\|_{L^{2}(\rr)}^{2}
\|\n_{\ep}^{\alpha/2}\na u_{\ep}\|_{L^{2}(\rr)}^{2}\right)\\
&\le C,
\ea\ee where we have used $\textbf{1}_{\{\n_{\ep}\ge1\}}\n_{\ep}^{(2-\alpha)/2}\le \textbf{1}_{\{\n_{\ep}\ge1\}}\n_{\ep}^{1/2}$ since $\alpha\ge1.$ The last two inequalities guarantees  \bnn \na m_{\ep} \in L^{2}\left(0,T;L^{1}(\rr)\right).\enn
Furthermore,  if   \be\la{y7}\p_{t}m_{\ep}\in L^{2}\left(0,T;W_{loc}^{-1,1}(\rr)\right).\ee
Then,  the Aubin-Lions lemma shows   there is a   $m\in L^{2}\left(0,T;L^{3/2}_{loc}(\rr)\right)$ such that  \be\la{s67}m_{\ep}\rightarrow m,\quad almost \,\,everywhere.\ee
This combining with \eqref{s53}   and    $\sqrt{\n_{\ep}}u_{\ep}\in L^{\infty}\left(0,T;L^{2}(\rr)\right)$  provides
\bnn \int_{\{\n_{\ep}\le 1\}} \frac{m^{2}}{\n^{2\alpha-1}}=\int_{\{\n_{\ep}\le 1\}} \liminf_{\ep\rightarrow0}\frac{(\n_{\ep}^{\alpha}u_{\ep})^{2}}{\n_{\ep}^{2\alpha-1}}\le \liminf_{\ep\rightarrow0}\int_{\rr}\n_{\ep}|u_{\ep}|^{2}\le C. \enn
So,   $m=0$ on vacuum sets. We define  $u=m\left(\chi(\n)\n^{\alpha}+(1-\chi(\n))\n^{(1+\alpha)/2}\right)^{-1}$ if $\n>0$ and $u=0$ if $\n=0$. The proof is thus completed.

We  need to check  \eqref{y7}. Let us first prove   \be\la{y18} \p_{t}\left((1-\chi(\n_{\ep}))\n_{\ep}^{(1+\alpha)/2}u_{\ep}\right)\in L^{2}\left(0,T;W^{-1,1}(\rr)\right).\ee
By \eqref{1} and \eqref{y13}, a careful   calculation shows
\be\la{y10}\ba &\p_{t}\left((1-\chi(\n_{\ep}))\n_{\ep}^{(1+\alpha)/2}u_{\ep}\right)\\
&=-\frac{2}{\alpha+1}\chi'\n_{\ep}u_{\ep}\p_{t}\n_{\ep}^{(\alpha+1)/2}  +(1-\chi(\n_{\ep}))\p_{t}(\n_{\ep}^{(1+\alpha)/2}u_{\ep})\\
&=\frac{2}{\alpha+1}\chi'\n_{\ep}u_{\ep}\left(\frac{\alpha-1}{2}\n_{\ep}^{(\alpha+1)/2} {\rm div} u_{\ep}+{\rm div} (\n_{\ep}^{(\alpha+1)/2}u_{\ep})\right)\\
 &\quad+(1-\chi(\n_{\ep}))\n_{\ep} u_{\ep}\left(\frac{3-\alpha}{2}\n_{\ep}^{(\alpha-1)/2} {\rm div} u_{\ep}-{\rm div} (\n_{\ep}^{(\alpha-1)/2}u_{\ep})\right)\\
&\quad +(1-\chi(\n_{\ep}))\n_{\ep}^{(\alpha-1)/2} \left[{\rm div}(\n_{\ep}^{\alpha}\na u_{\ep} + (\alpha-1)\n_{\ep}^{\alpha} {\rm div} u_{\ep}\mathbb{I})-\ep\na \n_{\ep}^{\ga}-{\rm div}(\n_{\ep}u_{\ep} \otimes u_{\ep})\right].\ea\ee
The terms in \eqref{y10} are dealt with as follows:   firstly,
\bnn\ba &\chi'\n_{\ep}u_{\ep}\left(\frac{\alpha-1}{2}\n_{\ep}^{(\alpha+1)/2} {\rm div} u_{\ep}+{\rm div} (\n_{\ep}^{(\alpha+1)/2}u_{\ep})\right)\\
&=\frac{\alpha+1}{2}\chi'  \n_{\ep}^{(\alpha+3)/2}u_{\ep} {\rm div}u_{\ep} +\frac{\alpha+1}{\alpha+3}u_{\ep}u_{\ep}\cdot \na \chi(\n_{\ep}^{(\alpha+3)/2})\\
&=\frac{\alpha+1}{2}\chi'  \n_{\ep}^{(\alpha+3)/2}u_{\ep} {\rm div}u_{\ep} +\frac{\alpha+1}{\alpha+3} (1-\chi(\n_{\ep}^{(\alpha+3)/2}))(u_{\ep}\cdot \na u_{\ep} +u_{\ep}{\rm div}u_{\ep})\\
&\quad -\frac{\alpha+1}{\alpha+3} \p_{j}((1-\chi(\n_{\ep}^{(\alpha+3)/2})) u_{\ep}^{k}  u_{\ep}^{j})\\
 &\in L^{2}\left(0,T;W^{-1,1}(\rr)\right),\ea\enn
where we have used \bnn\ba&\|(1-\chi(\n_{\ep}^{(\alpha+3)/2})) u_{\ep}^{k}  u_{\ep}^{j}\|_{L^{1}(\rr)}^{2}\\
 &+\int_{0}^{T}\|\n_{\ep}^{(\alpha+3)/2}\chi'(\n_{\ep})u_{\ep} {\rm div}u_{\ep}+(1-\chi(\n_{\ep}^{(\alpha+3)/2}))(u_{\ep}\cdot \na u_{\ep} +u_{\ep}{\rm div}u_{\ep})\|_{L^{1}(\rr)}^{2}\\
 &\le \|| u_{\ep}|^{2}\|_{L^{1}(\{1\le \n_{\ep}^{(\alpha+3)/2}\})}^{2}\\
  &\quad+\int_{0}^{T}\left(\||u_{\ep}|| {\rm div}u_{\ep}|\|_{L^{1}(\{1\le \n_{\ep}\le 2\})}^{2}+\||u_{\ep}||\na u_{\ep}|\|_{L^{1}(\{1\le \n_{\ep}^{(\alpha+3)/2}\})}^{2}\right)\\
&\le C \| \n_{\ep}  |u_{\ep}|^{2}\|_{L^{1}(\rr)}^{2}+ \int_{0}^{T} \| \sqrt{\n_{\ep}}  u_{\ep}  \|_{L^{2}(\rr)}^{2}\| \n_{\ep}^{\alpha/2} \na u_{\ep}  \|_{L^{2}(\rr)}^{2} \\
&\le C, \ea\enn owes to \eqref{7}.

Secondly, by  virtue of   \eqref{7} and \eqref{15s},
\be\la{y11}\ba &(1-\chi(\n_{\ep}))\n_{\ep}^{(\alpha-1)/2} {\rm div}(\n_{\ep}^{\alpha}\na u_{\ep})\\
&= {\rm div}((1-\chi(\n_{\ep}))\n_{\ep}^{\alpha-1/2}\n_{\ep}^{\alpha/2}\na u_{\ep})-\frac{\alpha-1}{2\alpha-1}(1-\chi(\n_{\ep}))\n_{\ep}^{\alpha/2}\na u_{\ep}\na \n_{\ep}^{\alpha-1/2}\\
&\quad+\chi'\n_{\ep}^{(3\alpha-1)/2}\na \n_{\ep}\cdot \na u_{\ep}\\
&\in L^{2}(0,T;W^{-1,1}(\rr)).\ea\ee
By similar argument, we   receive
\bnn (1-\chi(\n_{\ep}))\n_{\ep}^{(\alpha-1)/2} {\rm div} ( \n_{\ep}^{\alpha} {\rm div} u_{\ep}\mathbb{I})\in L^{2}(0,T;W^{-1,1}(\rr)).\enn

Next,
\bnn\ba &(1-\chi(\n_{\ep}))\n_{\ep} u_{\ep}\left(\frac{3-\alpha}{2}\n_{\ep}^{(\alpha-1)/2} {\rm div} u_{\ep}-{\rm div} (\n_{\ep}^{(\alpha-1)/2}u_{\ep})\right)\\
&\quad-(1-\chi(\n_{\ep}))\n_{\ep}^{(\alpha-1)/2}  {\rm div}(\n_{\ep}u_{\ep} \otimes u_{\ep}) \\
&=\frac{1-\alpha}{2} (1-\chi(\n_{\ep}))\n_{\ep}^{(\alpha+1)/2} u_{\ep} {\rm div} u_{\ep} -\frac{2}{\alpha+3}(1-\chi(\n_{\ep}^{(\alpha+3)/2}))(u_{\ep}{\rm div}u_{\ep}+u_{\ep}\cdot \na u_{\ep}) \\
&\quad+\p_{j}\left(\frac{2}{\alpha+3}(1-\chi(\n_{\ep}^{(\alpha+3)/2}))u_{\ep}^{k}u_{\ep}^{j} -(1-\chi(\n_{\ep}))u_{\ep}^{k}u_{\ep}^{j}  \n_{\ep}^{(\alpha+1)/2})\right)\\
&\in  L^{2}(0,T;W^{-1,1}(\rr)),\ea\enn
where the following inequality  has been used   \bnn \ba &\int_{0}^{T}\|(1-\chi(\n_{\ep}))\n_{\ep}^{(\alpha+1)/2} |u_{\ep}|^{2}\|_{L^{1}(\rr)}\\
&\le C\begin{cases}\int_{0}^{T}\|\na ((1-\chi(\n_{\ep}))\n_{\ep}^{(\alpha+1)/2} |u_{\ep}|)\|_{L^{1}(\rr)}^{2}+\|u_{\ep}\|_{L^{2}(\{\n_{\ep}\ge1\})}^{2}, & \,\,N=2\\
\int_{0}^{T}\|\n_{\ep}^{1/2} |u_{\ep}|^{2}\|_{L^{2}}^{2}+\|\n_{\ep}^{\alpha/2}\|_{L^{2}}^{2}, & \,\,N=3 \\\end{cases}\\
&\le C,\ea\enn  dues to  \eqref{e4}, \eqref{15s} and  \eqref{y17}.

Finally, since  \eqref{15s} and   $(2\ga+\alpha-1)/2\in (1,6\alpha-3)$, it has  \bnn\ba &(1-\chi(\n_{\ep}))\n_{\ep}^{(\alpha-1)/2}\na \n_{\ep}^{\ga}\\
&=\frac{2\ga}{2\ga+\alpha-1}\left(\na((1-\chi(\n_{\ep}))
\n_{\ep}^{(2\ga+\alpha-1)/2})+\n_{\ep}^{(2\ga+\alpha-1)/2}\chi\na \n_{\ep}\right)\\
&\in   L^{\infty}(0,T;W^{-1,1}(\rr)).\ea\enn

A similar argument yields
\bnn \p_{t}\left(\chi(\n_{\ep})\n_{\ep}^{\alpha}u_{\ep}\right)\in L^{2}\left(0,T;W^{-1,1}(\rr)\right),\enn
which combining with  \eqref{y18} gives  the desired  \eqref{y7}.

\emph{Step 2.} It satisfies
\be\la{s70} \sqrt{\n} u  \ln^{1/2} (e+|u|^{2})\in L^{\infty}\left(0,T;L^{2}(\rr)\right).\ee
To this end, let us first  check
\be\la{s62} \sqrt{\n_{\ep}} u_{\ep} \ln^{1/2} (e+|u_{\ep}|^{2})\in L^{\infty}\left(0,T;L^{2}(\rr)\right).\ee
  Clearly, the   \eqref{s62}  follows directly from \eqref{e4} and \eqref{7}  in case of $N=3.$  Now let us pay attention to   $N=2.$   Following   in \cite{mv,has1}, we have for any $\de\in (0,2)$
\be\la{s63}\ba &\frac{d}{dt}\int_{\rr}\n_{\ep}(1+|u_{\ep}|^{2})\ln(1+|u_{\ep}|^{2})
 + \int_{\rr}\n_{\ep}[1+\ln(1+|u_{\ep}|^{2})]|\na u_{\ep}|^{2}\\
 &\le C\int_{\rr}\n_{\ep}^{\alpha} |\na u_{\ep}|^{2}+C\ep^{2}\left(\int_{\rr}\n_{\ep}^{(4\ga-2\alpha-\de)/(2-\de)}\right)^{\frac{2-\de}{2}}
 \left(\int_{\rr}\n_{\ep}(1+|u_{\ep}|^{2})\right)^{\de/2}.\ea\ee
Thus, using \eqref{70} and \eqref{7},  integration of   \eqref{s63} in time conclude  the \eqref{s62}, so long as
 \bnn\ep^{2}\int_{0}^{T}\left(\int_{\rr}\n_{\ep}^{(4\ga-2\alpha-\de)/(2-\de)}\right)^{\frac{2-\de}{2}}\le C,\enn
which is fulfilled  because of \eqref{15s}.
Making use of \eqref{s53}, \eqref{y20}, \eqref{s62},  the Fatou Lemma, we  get  \eqref{s70}.

 \emph{Step 3.} Given  constant $M>1$, the \eqref{s53} and \eqref{y20} ensure that $\sqrt{\n_{\ep}}u_{\ep}|_{u_{\ep}\le M}\rightarrow \sqrt{\n}u|_{u\le M}$ almost everywhere   when  $\n>0$. If we also   define  $\sqrt{\n}u|_{u\le M}$ on sets   $\{\n=0\}$, then \bnn\sqrt{\n_{\ep}}u_{\ep}|_{u_{\ep}\le M}\le M\sqrt{\n_{\ep}}\rightarrow 0=\sqrt{\n}u|_{u\le M}\quad {\rm for}\quad \n=0.\enn
Recalling \eqref{15s}, it satisfies for  $q>2$
\bnn\sqrt{\n_{\ep}}u_{\ep}|_{u_{\ep}\le M}\in L^{\infty}\left(0,T;L^{q}\right),\enn
and therefore,  \bnn \int_{0}^{T}\int_{\rr}|\sqrt{\n_{\ep}}u_{\ep}|_{u_{\ep}\le M}-\sqrt{\n}u|_{u\le M}|^{2}\rightarrow 0\quad {\rm as}\,\,\,\,\ep\rightarrow 0.\enn
On the other hand, it follows from  \eqref{s62} and \eqref{s70} that
\bnn \ba&\int_{0}^{T}\int_{\rr}\left( \sqrt{\n_{\ep}}u_{\ep}|_{u_{\ep}> M}+\sqrt{\n}u|_{u> M} \right)^{2}\\
&\le \frac{C}{\ln(1+M^{2})}\int_{0}^{T}\int_{\rr} \left(\n_{\ep} |u_{\ep}|^{2}\ln(1+|u_{\ep}|^{2})+\n|u|^{2}\ln(1+|u|^{2})\right)\\
&\rightarrow 0\quad {\rm as}\,\,\,\,M\rightarrow \infty.\ea\enn
In conclusion, sending   $\ep\rightarrow 0$ first  and then  $M\rightarrow \infty$  yields  \bnn\ba \int_{0}^{T}\int_{\rr}|\sqrt{\n_{\ep}}u_{\ep}-\sqrt{\n}u|^{2}
&\le 2\int_{0}^{T}\int_{\rr}|\sqrt{\n_{\ep}}u_{\ep}|_{u_{\ep}\le M}-\sqrt{\n}u|_{u\le M}|^{2} \\
&\quad +2\int_{0}^{T}\int_{\rr}\left(|\sqrt{\n_{\ep}}u_{\ep}|_{u_{\ep}> M}|^{2}+|\sqrt{\n}u|_{u> M}|^{2}\right)\rightarrow 0.\ea\enn

\begin{lemma} \la{lem3.2} It satisfies
\be\la{pp0} \n_{\ep} u_{\ep}\rightarrow  \n  u \quad  in\quad L^{2}\left(0,T;L_{loc}^{q_{2}}(\rr)\right),\ee
where $q_{2}\in [1,2)$ if $N=2$ and $q_{2}\in [1,\frac{12\alpha-6}{6\alpha-1})$ if $N=3.$ \end{lemma}
\emph{Proof.}
Making use of  Lemma \ref{lem3.3},
 \eqref{15s},  \eqref{s70},  and the  inequality $$|\sqrt{x}-\sqrt{y}|\le \sqrt{|x-y|},\quad\forall\quad x\ge0,\,\,y\ge0,$$ we conclude the  \eqref{pp0} from the following
\be\la{pp}\ba &\|\n_{\ep}u_{\ep}-\n u\|_{L^{q_{2}}(\rr)}\\
&\le \|\sqrt{\n_{\ep}}(\sqrt{\n_{\ep}}u_{\ep}-\sqrt{\n} u)\|_{L^{q_{2}}(\rr)}+\|(\sqrt{\n_{\ep}}-\sqrt{\n})\sqrt{\n}u\|_{L^{q_{2}}(\rr)}\\
&\le \|\sqrt{\n_{\ep}}\|_{L^{q}(\rr)}\|(\sqrt{\n_{\ep}}u_{\ep}-\sqrt{\n} u)\|_{L^{2}(\rr)}+\|\sqrt{\n_{\ep}}-\sqrt{\n}\|_{L^{q}(\rr)}\|\sqrt{\n}u\|_{L^{2}
(\rr)}\\
&\le C\left(\|(\sqrt{\n_{\ep}}u_{\ep}-\sqrt{\n} u)\|_{L^{2}(\rr)}+\|\n_{\ep}-\n\|_{L^{q/2}(\rr)}^{1/2}\right),
\ea\ee
where  $q_{2}=(1/q+1/2)^{-1}$ with  $q\ge2$ if $N=2$ and  $q\in [2,6\alpha-3)$ if $N=3.$

\section{Proof of Theorem \ref{t2}}
Utilizing \eqref{6}, we deduce from   \eqref{1a}  that
\be\la{70x} \|\n(\cdot,t)\|_{L^{1}(\rr)}=\|\n_{\ep}(\cdot,t)\|_{L^{1}(\rr)}=\|\n_{0}\|_{L^{1}(\rr)},\ee
\be\la{7x}\ba&\sup_{t\in [0,T]}\int_{\rr}\n|u|^{2} +\int^{T}_{0}\int_{\rr} \left(\n ^{\alpha}|\na u |^{2}+(\alpha-1)\n ^{\alpha}({\rm div} u)^{2}\right)\le \int_{\rr}\frac{|m_{0}|^{2}}{\n_{0}}\ea\ee
and
\be\la{8x}\ba \sup_{t\in [0,T]}\int_{\rr} \n |u +\alpha \n ^{\alpha-2}\na \n |^{2}\le \int_{\rr}\n_{0}|u_{0}+\alpha \n_{0}^{\alpha-2}\na \n_{0}|^{2}.\ea\ee
In addition,   the same method  as  \eqref{15s} runs
 \be\ba \la{15ss} &\n \in L^{\infty}\left(0,T; L^{1}(\rr)\cap L^{q}(\rr) \right),\quad
 \na \n^{\alpha-1/2}\in L^{\infty}\left(0,T; L^{2}(\rr)\right),\ea \ee
where    $ q<\infty $ if $N=2 $ and $q=6\alpha-3$ if $N=3$.

Set $z=\n_{\ep}-\n$. Subtracting  $\eqref{1a}_{1}$ from  $\eqref{1}_{1}$  receives
\bnn\begin{cases}
\p_{t}z=\triangle\left ( \n_{\ep}^{\alpha}- \n^{\alpha}\right)-{\rm div} \left[\n_{\ep}(u_{\ep}+\alpha \n_{\ep}^{\alpha-2}\na \n_{\ep})-\n(u+\alpha \n^{\alpha-2}\na \n)\right],\\
z(x,t=0)=0.\end{cases}\enn
Multiplying the above  by $\vp(x,t)\in C_{0}^{\infty}(\rr\times [0,+\infty))$ and integrating the expression by parts give rise to
\be\ba\la{10}  \int_{\rr}z \vp(x,t)
&=\int_{0}^{t}\int_{\rr} z\left(\vp_{s}+a\triangle \vp\right)\\
&\quad+\int_{0}^{t}\int_{\rr}\left[\n_{\ep}(u_{\ep}+\alpha \n_{\ep}^{\alpha-2}\na \n_{\ep})-\n(u+\alpha \n^{\alpha-2}\na \n)\right]\cdot \na \vp,\ea\ee
where $ a=(\n_{\ep}^{\alpha}-\n^{\alpha})/z$ if $z\neq 0$  and $a=
0$ if $z=0.$

To be continued,    consider the  following  backward   parabolic equation
\be\la{13}\begin{cases}
\p_{s} \vp_{R}+a_{n}\triangle \vp_{R}=0,\quad& x\in B_{R},\,\,0\le s<t,\\
\vp_{R}=0,\quad &x\in \p B_{R},\,\ 0\le s<t,\\
\vp_{R}(x,t)=\th(x)\in H_{0}^{1}(B_{R}),& x\in  B_{R},
\end{cases}\ee
where   $ a_{n}=\eta_{1/n}*a_{K,\ep}\in [2^{-1}\ep,2K]$  and $\eta_{1/n}$ being  the standard Friedrichs' mollifier such that  as $n\rightarrow \infty,$ \be\la{18}a_{n}\rightarrow  a_{K,\ep}= \begin{cases}K,& a>K,\\
a,& \ep\le a\le K,\\
\ep, &a<\ep. \end{cases} \ee The classical  linear parabolic theory (cf. \cite{la}) ensures  that   \eqref{13} has a unique   solution $\vp_{R}\in L^{\infty}\left(0,t; H^{1}_{0}\right)\cap L^{2}\left(0,t; H^{2}\right)$.  If we
multiply  \eqref{13} by $\triangle \vp_{R}$, we infer    for  $\tau \in [0,t]$ \bnn\ba \frac{1}{2}\int_{B_{R}}  |\na \vp_{R}|^{2}(x,\tau) &-\frac{1}{2}\int_{B_{R}}  |\na \vp_{R}|^{2}(x,t) +\int_{\tau}^{t}\int_{B_{R}} a_{n} |\triangle \vp_{R}|^{2}=0,\ea\enn
and thus,
\be\la{75}\ba  \int_{B_{R}}  |\na \vp_{R}|^{2}(x,\tau)+\int_{\tau}^{t}\int_{B_{R}} a_{n} |\triangle \vp_{R}|^{2}\le \int_{B_{R}}  |\na \th|^{2}. \ea\ee

Define a smooth cut-off function $\xi_{R}$ satisfying
\be\la{19} \xi_{R}=1\,\,\,{\rm in }\,\,B_{R/2},\quad \xi_{R}=0 \,\,\,{\rm in }\,\,\rr\setminus B_{R},\quad   |\na^{k} \xi_{R}|\le CR^{-k},\,(k=1,2).\ee
 If we  extend  $\vp_{R}$ to  $\rr$ by  zero and  replace $\vp_{R}$ in \eqref{10} with   $\vp=\xi_{R}\vp_{R}$,  the first term on the right-hand side of \eqref{10} satisfies
\be\la{29a}\ba &\int_{0}^{t}\int_{\rr} z\left(\vp_{s}+a\triangle \vp\right)\\
&=\int_{0}^{t}\int_{\rr} z\xi_{R} \left(\p_{s}\vp_{R}+a\triangle \vp_{R}\right)+\int_{0}^{t}\int_{\rr} (\n_{\ep}^{\alpha}-\n^{\alpha}) \left(2\na \xi_{R}\na\vp_{R}+\vp_{R}\triangle \xi_{R}\right)\\
&=\int_{0}^{t}\int_{\rr} z\xi_{R} \left(a-a_{n}\right)\triangle \vp_{R}+\int_{0}^{t}\int_{\rr} (\n_{\ep}^{\alpha}-\n^{\alpha}) \left(2\na \xi_{R}\na\vp_{R}+\vp_{R}\triangle \xi_{R}\right)\\
&\triangleq I_{1} +I_{2},\ea\ee where in the second equality we used   \eqref{13}.

By \eqref{18} and \eqref{75},  one has
\bnn\ba
|I_{1}|&\le \left(\int_{0}^{t}\int_{B_{R}} \frac{z^{2}\left(a-a_{n}\right)^{2}}{a_{n}}\right)^{1/2} \left(\int_{0}^{t}\int_{B_{R}} a_{n}|\triangle \vp_{R}|^{2}\right)^{1/2}\\
&\le  \|\na \th\|_{L^{2}(B_{R})}\left(\int_{0}^{t}\int_{B_{R}} \frac{z^{2}\left(a-a_{n}\right)^{2}}{a_{n}}\right)^{1/2}\\
&\le \sqrt{2\ep^{-1}}\|\na \th\|_{L^{2}(B_{R})} \left(\int_{0}^{t}\int_{B_{R}}  z^{2}\left(a-a_{K,\ep}\right)^{2}+ z^{2}\left(a_{K,\ep}-a_{n}\right)^{2}\right)^{1/2}\\
&\le  C \ep^{1/2}\|\na \th\|_{L^{2}(B_{R})},\ea\enn
where the last inequality owes to \eqref{15s}, \eqref{15ss},   and the following two inequalities:
 \bnn\ba \int_{0}^{t}\int_{B_{R}}  z^{2}\left(a_{K,\ep}-a_{n}\right)^{2}\le C\int_{0}^{T}\|z\|_{L^{6\alpha-3}}^{2} \|(a_{K,\ep}-a_{n})\|_{L^{\frac{12\alpha-6}{6\alpha-5}}}^{2}\rightarrow0 \quad(n\rightarrow\infty)\ea\enn
 and
\bnn\ba  \int_{0}^{t}\int_{B_{R}}z^{2}\left(a-a_{K,\ep}\right)^{2}
 &\le C\ep^{2} \sup_{0\le t\le T}\|z\|_{L^{2}}^{2}+\int_{0}^{T}\int_{B_{R}\cap \{x\in \rr: a>K\}}   \left(\n_{\ep}^{\alpha}-\n^{\alpha}\right)^{2}\\
  &\le C\ep^{2} \quad(K\rightarrow\infty).\ea\enn
Therefore,  \eqref{29a} is estimated  as
\be\la{29x}\ba \int_{0}^{t}\int_{\rr} z\left(\vp_{s}+a\triangle \vp\right) \le C \ep^{1/2}\|\na \th\|_{L^{2}(B_{R})} + I_{2}.\ea\ee

The second  term  on  the   right-hand side of \eqref{10} satisfies for  $\vp=\xi_{R}\vp_{R}$
\be\la{72}\ba &\int_{0}^{t}\int_{\rr}\left[\n_{\ep}(u_{\ep}+\alpha \n_{\ep}^{\alpha-2}\na \n_{\ep})-\n(u+\alpha \n^{\alpha-2}\na \n)\right]\cdot \na (\xi_{R}\vp_{R})\\
&\le C\int_{0}^{T}\int_{B_{R}} \left(\n_{\ep}|u_{\ep}+\alpha \n_{\ep}^{\alpha-2}\na \n_{\ep}|+\n|u+\alpha \n^{\alpha-2}\na \n|\right) |\na \vp_{R}|\\
&\quad +\int_{0}^{T}\int_{B_{R}}\left(\n_{\ep}|u_{\ep}+\alpha \n_{\ep}^{\alpha-2}\na \n_{\ep}|+\n|u+\alpha \n^{\alpha-2}\na \n|\right)|\na \xi_{R}||\vp_{R}|\\
&\triangleq J_{1}+J_{2}.\ea\ee
Owing to  \eqref{11} and \eqref{8},
\be\ba\la{12} \|\sqrt{\n_{\ep}}|u_{\ep}+\alpha \n_{\ep}^{\alpha-2}\na \n_{\ep}|\|_{L^{2}(\rr)}^{2}\le  \ep \|\n_{0}\|_{L^{\ga}(\rr)}^{\ga}\le C \ep.\ea\ee
This together  with \eqref{11} and \eqref{8x}   shows for   $\bar{q}>2$
 \be\ba\la{16} J_{1}&\le  C\int_{0}^{T}\|\sqrt{\n_{\ep}}|u_{\ep}+\alpha \n_{\ep}^{\alpha-2}\na \n_{\ep}|\|_{L^{2}(B_{R})}\|\sqrt{\n_{\ep}}\|_{L^{(1/2-1/\bar{q})^{-1}}(B_{R})}\|\na \vp_{R}\|_{L^{\bar{q}}(B_{R})}\\
&\le  C\ep^{1/2}\int_{0}^{T}\|\sqrt{\n_{\ep}}\|_{L^{(1/2-1/\bar{q})^{-1}}(B_{R})}\|\na \vp_{R}\|_{L^{\bar{q}}(B_{R})}.
 \ea\ee

 We discuss $J_{1}$ in two cases.
\begin{itemize}
  \item \underline{ \emph{Let $\bar{q}=2+\de$ with $\delta>0$ small  in case of $N=2$.}}
\end{itemize}

By \eqref{15s} and  \eqref{lpo},    \be\la{17}\ba  J_{1}
&\le  C\ep^{1/2}\int_{0}^{T}\|\sqrt{\n_{\ep}}\|_{L^{(1/2-1/(2+\de))^{-1}}(B_{R})}\|\na \vp_{R}\|_{L^{2+\de}(B_{R})} \\
&\le  C\ep^{1/2}\int_{0}^{T} \|\na \vp_{R}\|_{L^{2+\de}(B_{R})}\\
&\le C\ep^{1/2}  \int_{0}^{T}\left(\|\na \vp_{R}\|_{L^{2}(B_{R})}+\ep^{\frac{-\delta}{2(2+\delta)}}\|\na \vp_{R}\|_{L^{2}(B_{R})}^{\frac{2}{2+\delta}}\|\sqrt{a_{n}}\triangle  \vp_{R}\|_{L^{2}(B_{R})}^{\frac{\delta}{2+\delta}}\right)\\
&\le C\ep^{\frac{1}{2+\delta}}\left(\sup_{t\in [0,T]}\|\na \vp_{R}\|_{L^{2}(B_{R})}+\left(\int_{0}^{T}\|\sqrt{a_{n}}\triangle \vp_{R}\|_{L^{2}(B_{R})}^{2}\right)^{1/2}\right)\\
 &\le C\ep^{\frac{1}{2+\delta}}\|\na \th\|_{L^{2}(B_{R})},
\ea\ee
where  in the   third inequality  we  have   used
 \be\la{y41}\|\na^{2} \vp_{R}\|_{L^{2}(B_{R})}\le C\|\triangle \vp_{R}\|_{L^{2}(B_{R})}\le \ep^{-1/2}\|\sqrt{a_{n}}\triangle \vp_{R}\|_{L^{2}(B_{R})},\ee  owes to \eqref{18} and Lemma \ref{lem2.2}.

\begin{itemize}
  \item \underline{\emph{Let $\bar{q}=\frac{6\alpha-3}{3\alpha-2}$ in  case of $N=3$.}}
\end{itemize}

   Similar to \eqref{17}, we deduce
 \be\la{35}\ba  J_{1}
&\le  C\ep^{1/2}\int_{0}^{T}\|\sqrt{\n_{\ep}}\|_{L^{12\alpha-6}(B_{R})}\|\na \vp_{R}\|_{L^{\frac{6\alpha-3}{3\alpha-2}}(B_{R})}\\
&\le C\ep^{1/2}  \int_{0}^{T}\left((1+\ep^{\frac{-1}{4(2\alpha-1)}})\|\na \vp_{R}\|_{L^{2}(B_{R})}+\ep^{\frac{-1}{4(2\alpha-1)}}\|\sqrt{a_{n}}\triangle \vp_{R}\|_{L^{2}(B_{R})}\right)\\
&\le C\ep^{\frac{4\alpha-3}{4(2\alpha-1)}}\|\na \th\|_{L^{2}(B_{R})}.
\ea\ee

With the aid of   \eqref{16}, \eqref{17}, and \eqref{35},
the  \eqref{72} satisfies
\bnn\ba &\int_{0}^{t}\int_{\rr}\left(\n_{\ep}(u_{\ep}+\alpha \n_{\ep}^{\alpha-2}\na \n_{\ep})+\n(u+\alpha \n^{\alpha-2}\na \n)\right)\cdot \na (\xi_{R}\vp_{R})\\
&\le J_{2}+C\|\na \th\|_{L^{2}(B_{R})}\begin{cases}\ep^{\frac{1}{2+\delta}},&N=2, \\
\ep^{\frac{4\alpha-3}{4(2\alpha-1)}},&N=3,\end{cases}\ea\enn
which, along  with    \eqref{10}  and \eqref{29x}, implies
\be\ba\la{10x}  \int_{\rr}z\xi_{R}\vp_{R}(x,t) \le I_{2}+J_{2}+C\|\na \th\|_{L^{2}(B_{R})}\begin{cases}\ep^{\frac{1}{2+\delta}},&N=2, \\
\ep^{\frac{4\alpha-3}{4(2\alpha-1)}},&N=3.\end{cases}\ea\ee

Next, by  \eqref{15s}, \eqref{15ss},   \eqref{19},  the Poincar$\acute{e}$ inequality, we deduce from   \eqref{29a} that
\be\ba\la{xo} I_{2}
&=
\int_{0}^{t}\int_{\rr} (\n_{\ep}^{\alpha}-\n^{\alpha}) \left(2\na \xi_{R}\na\vp_{R}+\vp_{R}\triangle \xi_{R}\right)\\
&\le   CR^{-1}\int_{0}^{T}  \|\n_{\ep}^{\alpha}-\n^{\alpha}\|_{L^{2}(B_{R}  \backslash B_{R/2})} \| \na \vp_{R}\|_{L^{2}(B_{R})} \\
&\le   CR^{-1}  \| \na \th\|_{L^{2}(B_{R} )}\rightarrow0\quad (R\rightarrow\infty).\ea\ee
By  \eqref{11}, \eqref{8x}, \eqref{15s}, \eqref{75}, \eqref{19}, \eqref{12}, \eqref{y41}, the Poincar$\acute{e}$ inequality, we deduce from  \eqref{72} that
\be\la{xox}\ba  J_{2}
&\le  C\int_{0}^{T}\int_{B_{R}  \backslash B_{R/2}}\left(\n_{\ep}|u_{\ep}+\alpha \n_{\ep}^{\alpha-2}\na \n_{\ep}|+\n|u+\alpha \n^{\alpha-2}\na \n|\right)|\na \xi_{R} ||\vp_{R}|\\
&\le   CR^{-1}\int_{0}^{T} \|\sqrt{\n_{\ep}} \|_{L^{4}(B_{R}  \backslash B_{R/2})} \| \sqrt{\n_{\ep}}|u_{\ep}+\alpha \n_{\ep}^{\alpha-2}\na \n_{\ep}|\|_{L^{2}(\rr)} \| \vp_{R}\|_{L^{4}(B_{R})}\\
&\le   C \ep^{1/2}\sup_{t\in [0,T]} \| \n_{\ep} \|_{L^{2}(B_{R}  \backslash B_{R/2})}^{1/2} \int_{0}^{T}\| \na \vp_{R}\|_{L^{4}(B_{R})}\\
&\le   C\sup_{t\in [0,T]} \| \n_{\ep} \|_{L^{2}(B_{R}  \backslash B_{R/2})}^{1/2} \int_{0}^{T}  \left( \| \na \vp_{R}\|_{L^{2}(B_{R})}+\| \sqrt{a_{n}} \triangle \vp_{R}\|_{L^{2}(B_{R})}\right)\\
&\le   C\sup_{t\in [0,T]} \| \n_{\ep} \|_{L^{2}(B_{R}  \backslash B_{R/2})}^{1/2} \| \na \th\|_{L^{2}(B_{R} )}\rightarrow0\quad (R\rightarrow\infty).\ea\ee

Particularly, if we replace the function  $\vp_{R}(x,t)$ in \eqref{10x} with
 \be\la{r2} \vp_{R}(x,t)=\xi_{R}\left(\n_{\ep}^{\alpha-\frac{1}{2}}-\n^{\alpha-\frac{1}{2}}\right)(x,t),\ee
  we conclude  from  \eqref{10x}-\eqref{xox}   that by sending $R\rightarrow\infty$,
\bnn \int_{\rr}(\n_{\ep}-\n)(\n_{\ep}^{\alpha-\frac{1}{2}}-\n^{\alpha-\frac{1}{2}}) \le  C \begin{cases}\ep^{\frac{1}{2+\delta}},&N=2, \\
\ep^{\frac{4\alpha-3}{4(2\alpha-1)}},&N=3.\end{cases} \enn
In terms of  Lemma
 \ref{lem2.4}, it satisfies  for  $\alpha\ge \frac{3}{2}$ \bnn\int_{\rr}|\n_{\ep}-\n|^{\alpha+\frac{1}{2}}
 \le  C \begin{cases}\ep^{\frac{1}{2+\delta}},&N=2, \\
\ep^{\frac{4\alpha-3}{4(2\alpha-1)}},&N=3.\end{cases}\enn

The  proof of Theorem \ref{t2} is complete   by exploiting  \eqref{15s},  \eqref{15ss},  and  interpolation inequalities.

\section*{Appendix. Proof of \eqref{e4}}
   Multiplying   equations $\eqref{1}_{2}$ by $4|u_{\ep}|^{2}u_{\ep}$ yields
\bnn\ba &\frac{d}{dt}\int_{\mathbb{R}^{3}}\n_{\ep}|u_{\ep}|^{4}+4\int_{\mathbb{R}^{3}}\n^{\alpha}_{\ep}|u_{\ep}|^{2}|\na u_{\ep}|^{2}\\
&\quad +8\int_{\mathbb{R}^{3}}\n^{\alpha}_{\ep}|u_{\ep}|^{2}|\na |u_{\ep}| |^{2}+4(\alpha-1)\int_{\mathbb{R}^{3}}\n^{\alpha}_{\ep}{\rm div}u_{\ep} |u_{\ep}| \left({\rm div}u_{\ep}|u_{\ep}| +2 u_{\ep}\cdot \na |u_{\ep}|\right)\\
&=4\ep\int_{\mathbb{R}^{3}}\na \n_{\ep}^{\ga}{\rm div}(|u_{\ep}|^{2}u_{\ep}).\ea\enn
By Young's inequality, it satisfies for all $\alpha>1/2$
\bnn\ba &8\int_{\mathbb{R}^{3}}\n^{\alpha}_{\ep}|u_{\ep}|^{2}|\na |u_{\ep}| |^{2}+4(\alpha-1)\int_{\mathbb{R}^{3}}\n^{\alpha}_{\ep}{\rm div}u_{\ep} |u_{\ep}| \left({\rm div}u_{\ep}|u_{\ep}| +2 u_{\ep}\cdot \na |u_{\ep}|\right)\\
&\ge-\frac{5}{2}\int_{\mathbb{R}^{3}} \n^{\alpha}_{\ep}({\rm div}u_{\ep})^{2} |u_{\ep}|^{2}\ge- \frac{5}{2}\int_{\mathbb{R}^{3}} \n^{\alpha}_{\ep}|\na u_{\ep}|^{2} |u_{\ep}|^{2}.\ea\enn
Next,
\bnn\ba &\ep\int_{\mathbb{R}^{3}}\na \n_{\ep}^{\ga}{\rm div}(|u_{\ep}|^{2}u_{\ep})\\
&\le  \ep \int_{\mathbb{R}^{3}}\n^{\alpha}_{\ep}|u_{\ep}|^{2}|\na u_{\ep}|^{2} +C\ep \int_{\mathbb{R}^{3}}\n_{\ep}^{2\ga-\alpha}|u_{\ep}|^{2}\\
&\le  \ep \int_{\mathbb{R}^{3}}\n^{\alpha}_{\ep}|u_{\ep}|^{2}|\na u_{\ep}|^{2}+ C\ep \|\n_{\ep}^{2\ga-\alpha-1/2}\|_{L^{2}(\mathbb{R}^{3})}
\left(1+\|\n_{\ep}^{1/2}|u_{\ep}|^{2}\|_{L^{2}(\mathbb{R}^{3})}^{2}\right).\ea\enn
For small $\ep\le 1/2,$ the above three inequalities ensure that
\bnn\ba &\frac{d}{dt}\int_{\mathbb{R}^{3}}\n_{\ep}|u_{\ep}|^{4}+ \int_{\mathbb{R}^{3}}\n^{\alpha}_{\ep}|u_{\ep}|^{2}|\na u_{\ep}|^{2}\\
&\le C\ep \|\n_{\ep}^{2\ga-\alpha-1/2}\|_{L^{2}(\mathbb{R}^{3})}
\left(1+\|\n_{\ep}^{1/2}|u_{\ep}|^{2}\|_{L^{2}(\mathbb{R}^{3})}^{2}\right).\ea\enn
The proof can be done by means of the Gronwall inequality, provided
\bnn \ep\int_{0}^{T}\|\n_{\ep}^{2\ga-\alpha-1/2}\|_{L^{2}(\mathbb{R}^{3})}dt\le C.\enn
In fact,
since \eqref{e1} and \eqref{e2} implies $1\le 4\ga-2\alpha-1\le 2\ga+4\alpha-3,$ it has
\bnn\ba \|\n_{\ep}^{2\ga-\alpha-1/2}\|_{L^{2}(\mathbb{R}^{3})}^{2}&=\left(\int_{\{\n_{\ep}\le1\}}
+\int_{\{\n_{\ep}\ge1\}}\right)\n_{\ep}^{4\ga-2\alpha-1}\\
&\le  \int_{\{\n_{\ep}\le1\}}\n_{\ep}+\int_{\{\n_{\ep}\ge1\}} \n_{\ep}^{2\ga+4\alpha-3} \le C+ C \|\na \n_{\ep}^{(\ga+\alpha-1)/2}\|_{L^{2}(\mathbb{R}^{3})}^{4}.\ea\enn
where the last inequality owes to \eqref{15s}, Sobolev inequality and the following
\bnn\ba \int_{\mathbb{R}^{3}} \n_{\ep}^{2\ga+4\alpha-3}
 \le C\|\n_{\ep}\|_{L^{6\alpha-3}(\mathbb{R}^{3})}^{2\alpha-1}\|\n_{\ep}^{(\ga+\alpha-1)/2}\|_{L^{6}(\mathbb{R}^{3})}^{4} \le C \|\na \n_{\ep}^{(\ga+\alpha-1)/2}\|_{L^{2}(\mathbb{R}^{3})}^{4}.\ea\enn
Therefore,
\bnn \ba  \ep\int_{0}^{T}\|\n_{\ep}^{2\ga-\alpha-1/2}\|_{L^{2}(\mathbb{R}^{3})} \le C\ep +\ep\int_{0}^{T}\|\na \n_{\ep}^{(\ga+\alpha-1)/2}\|_{L^{2}(\mathbb{R}^{3})}^{2} \le C,\ea\enn
where the last inequality owes to \eqref{8} and the $C$ is independent of $\ep$.

\end{document}